\documentclass{article}
\usepackage{amscd,amssymb,amsmath,verbatim,amsthm}
\usepackage{comment}

\newtheorem{thm}{Theorem}[section]
\newtheorem{cor}[thm]{Corollary}
\newtheorem{lem}[thm]{Lemma}
\newtheorem{prop}[thm]{Proposition}
\newtheorem{claim}[thm]{Claim}
\newtheorem{defn}[thm]{Definition}
\theoremstyle{remark}

\numberwithin{equation}{section}

\newcommand{\dist}{\mathrm{dist}}

\newcommand{\vol}{\mathrm{Vol}}

\newcommand{\inj}{\mathrm{Inj}}

\def\XXint#1#2#3{{\setbox0=\hbox{$#1{#2#3}{\int}$}
     \vcenter{\hbox{$#2#3$}}\kern-.5\wd0}}

\newcommand{\RR}{\mathbb R}

\newcommand{\del}{\partial}

\newcommand{\Ric}{\mathrm{Ric}}

\newcommand{\calC}{{\mathcal C}}

\newcommand{\calE}{{\mathcal E}}

\newcommand{\calN}{{\mathcal N}}
\newcommand{\calO}{{\mathcal O}}

\newcounter{mnotecount}[section] 

\date{}

\begin{document}

\title{Ricci flow on asymptotically conical surfaces with nontrivial topology}
\author{James Isenberg \thanks{Research supported in part by NSF grant PHY-0652903 } \\ University of Oregon \and Rafe Mazzeo 
\thanks{Research supported in part by NSF grant DMS- 0805529}\\ Stanford University
\and Natasa Sesum  \thanks{Research supported in part by NSF grant DMS-0905749} \\ University of Pennsylvania}

\maketitle

\begin{abstract}
As part of the general investigation of Ricci flow on complete surfaces with finite total curvature, we study this flow for surfaces
with asymptotically conical (which includes as a special case asymptotically Euclidean) geometries. After establishing long-time
existence, and in particular the fact that the flow preserves the asymptotically conic geometry, we prove that the solution metric
$g(t)$ expands at a locally uniform linear rate; moreover, the rescaled family of metrics $t^{-1}g(t)$ exhibits a transition
at infinite time inasmuch as it converges locally uniformly to a complete, finite area hyperbolic metric which is the unique
uniformizing metric in the conformal class of the initial metric $g_0$. 
\end{abstract}

\section{Introduction}
This paper is a continuation of an ongoing general investigation of the global properties of Ricci flow on various 
classes of complete surfaces.  Earlier work includes the paper \cite{JMS} by  L.\ Ji together with the second and third authors 
here, which studies Ricci flow surfaces with asymptotically hyperbolic cusp ends and negative Euler characteristic,
and the paper of Albin, Aldana and Rochon \cite{AAR} concerning the flow on surfaces with infinite area asymptotically
hyperbolic ends.  In the present paper we study the behaviour of this flow on surfaces with asymptotically Euclidean ends, 
or slightly more generally, asymptotically conic ends, and with $\chi(M) < 0$. 

An open complete Riemannian surface $(M,g)$ is called asymptotically conical (or AC for short) of order $\tau > 0$ 
if $M$ is topologically finite (i.e., has finite genus and a finite number of ends) and if on each end $E$ of $M$, $g$ 
is asymptotic to a model conic metric $dr^2 + \alpha^2 r^2 d\theta^2$ of angle $2\pi \alpha$, for some constant
$\alpha > 0$ with an error term which decays like $r^{-\tau}$ at infinity along with a certain number of its derivatives.  
We give a more precise definition in the next section.

The general problem studied here is to determine the long-time behaviour of the Ricci flow 
\begin{equation}
\del_t g = -2 \Ric(g(t)), \qquad g(0) = g_0,
\label{eq:rfs0}
\end{equation}
where the initial metric $g_0$ is asymptotically conical and $\dim M = 2$. Prior work in this setting includes the older 
paper by Wu \cite{Wu}  and the more recent work by Javaheri and the first author \cite{IJ}, in both of which the underlying 
surface is diffeomorphic to $\RR^2$. Wu's paper allows initial metrics $g_0$ with complicated asymptotic behaviour 
at infinity. In either case, the conclusion is that under appropriate hypotheses, the flow converges to either 
flat $\RR^2$ or else the `cigar soliton'. Here we allow general topologies, but impose more stringent conditions on the 
asymptotic geometry. Some local aspects of this flow for asymptotically Euclidean (or asymptotically locally Euclidean--ALE) 
metrics in higher dimensions have been obtained by Oliynyk and Woolgar \cite{OW} and by Dai and Ma \cite{DM}.  The first of 
these papers proves some results concerning long-time existence and convergence in the rotationally symmetric case, and the latter 
paper establishes, amongst other things, monotonicity of the ADM mass when $\dim M \geq 3$. Another interesting stability
result with weaker conditions imposed on the asymptotic behaviour of the initial metric 
has been proven by Schn\"urer, Schulze and Simon \cite{SSS}. 

For $\dim M = 2$, one expects to be able to obtain very explicit convergence results for Ricci flow. However, there are also
obstructions which indicate the possibility of some interesting new phenomena. For example, in the case considered
here, if $\chi(M) < 0$ and if $(M,g_0)$ is asymptotically conic, then the conformal type of $(M,g_0)$ is that of a punctured 
Riemann surface and the (unique) uniformizing metric is a complete hyperbolic metric with finite area (and hence cusp ends). 
Therefore if in this case the flow $(M, g(t))$ exists for all times and is asymptotically conic for each $t > 0$, and if it converges in certain sense then it seems
that the asymptotic geometry must change drastically at $t=\infty$. This intriguing scenario motivates our main result, which we now state:
\begin{thm}
Let $(M,g_0)$ be a surface with asymptotically conical ends and $\chi(M) < 0$. If  the Ricci evolution $g(t)$ is written in the form  $u(t)g_0$, 
then there are  constants $C_1, C_2 > 0$ depending only on $g_0$ such that $C_1 \le u(\cdot,t) \leq C_2(1+t)$ for all $t \geq 0$. 
In addition, for each compact set $K \subset \subset M$ there is a constant $C_K > 0$ such that $C_K (1+t) \leq u(t,x)$ 
for all $x \in K$ and $t \geq 0$.  

If one defines the rescaled metric $\tilde{g}(t) = t^{-1} g(t)$, then $\tilde{g}$ converges smoothly on each compact
set to a limiting metric $\tilde{g}_\infty$ on $M$ which is a complete hyperbolic metric with finite area, and hence is the unique
uniformizing metric in the conformal class of $g_0$. 
\end{thm} 

A similar rescaling by $1/t$ to normalize the metric and to obtain a smooth limit has been used earlier by Lott \cite{Lo}
in an analysis of a three-dimensional flow.  Recent work by Dai and Z.\ Zhang  (in preparation) studies Ricci flow in this 
same two-dimensional setting; they obtain pointed Gromov-Hausdorff convergence of the unrescaled metric to a flat 
plane or cigar soliton. 

The plan of this paper is as follows.  We begin by recalling some background information about AC geometry and the Ricci flow on 
surfaces. We then indicate a proof of short-time existence within the class of asymptotically conic metrics. This is fairly standard 
and has been proved elsewhere, so we sketch this quite briefly. Long-time existence of solutions to the flow is proved using the 
standard device of finding a potential function; this argument is similar to, but simpler than, the corresponding argument in \cite{JMS}. 
The linear in time upper bound and uniform lower bound for $u$ follow from simple barrier estimates, but the locally uniform linear in time
lower bound emerges only as a consequence of a fairly lengthy and quite geometric argument which brings in a number of
other ingredients. Once we have obtained this linear lower bound, the proof of convergence for the rescaled family of metrics 
$\tilde{g}(t)$ is straightforward.

\section{Preliminaries}
This section contains a few geometric facts about asymptotically conic surfaces and analytic facts about the Ricci flow which 
we use in our study.

Let us start with a more careful discussion of the space of asymptotically conical metrics on an open surface $M$.
For any number $\alpha > 0$, we define the model cone with angle $2\pi \alpha$ by
\[
C_\alpha = \RR^+ \times S^1, \qquad g_\alpha = dr^2 + \alpha^2 r^2 d\theta^2.
\]
Now, let $M$ be an open surface with finite topology, and suppose that we have chosen a fixed identification of each end $E_j$ 
of $M$ with the `open' end of some $C_\alpha$. In other words, suppose that we have fixed  coordinates $r \geq R$ and $\theta \in S^1$ on each $E_j$. 
A metric $g$ on $M$ is called asymptotically conical (AC) if for each end $E_j$,
\[
\left. g \right|_{E_j} = g_{\alpha_j} + k_j, \qquad \mbox{where}\qquad |k_j|_{g_{\alpha_j}} \leq C r^{-\tau}
\]
for some $\tau > 0$. More precisely, we suppose that $k_j$ lies in a weighted H\"older space, defined as follows.
Let $\Lambda^{k,\alpha}(M)$ denote the usual H\"older space of order $k + \alpha$ defined relative to some fixed
background metric $\bar{g}$ which is exactly equal to $g_{\alpha_j}$ on each $E_j$. This definition extends immediately
to sections of any Hermitian bundle $F$ over $M$ which has a fixed trivialization over each $E_j$ (e.g., \ a tensor bundle 
with the induced metric). Next, for any real number $\mu$, let $r^\mu \Lambda^{k,\alpha}(M,F)$ denote the space of sections 
$u$ of the form $u = r^\mu \tilde{u}$ where $\tilde{u} \in \Lambda^{k,\alpha}(M,F)$. Here and below, we fix a function $r$ on 
$M$ which is everywhere smooth and strictly positive and which agrees with the chosen radial function $r$ on each end. With the above notation we have the following definition.
\begin{defn}
The metric $g$ on $M$ is called asymptotically conical (or AC) of order $\tau$ if on each end $E_j$, the error term
$k_j$ lies in the space $r^{-\tau}\Lambda^{2,\alpha}(M, S^2 T^*M)$.
\end{defn}
It follows immediately that if $g$ is AC of order $\tau$, then its scalar curvature $R = R_g$ lies in the
space $r^{-2-\tau}\Lambda^{0,\alpha}(M)$. 

Shi (\cite{Sh}) has proven that if $(M,g_0)$ is a complete noncompact manifold with bounded curvature (which is true for our initial metric $g_0$) then the Ricci flow  (\ref{eq:rfs0}) has a solution with
bounded curvature on a short time interval. In \cite{CZ} it has been proven that a solution with a particular set of prescribed initial data  is unique within the class of solutions with  bounded curvatures in space. It follows  that we have  short time existence for a solution to (\ref{eq:rfs0})  starting with the asymptotically conical metric $g_0$ on $M$. Our goal is to understand the long time behavior of our solution $g(\cdot,t)$. 

The general Ricci flow equation \eqref{eq:rfs0} takes a particularly simple form in two dimensions. In this setting, 
$\Ric(g) = \frac12 R \, g$ where $R$ is the scalar curvature (i.e.\ twice the Gauss curvature), so \eqref{eq:rfs0}
is the same as 
\begin{equation}
\del_t g_{ij}(t) = - R g_{ij}(t);
\label{eq:rfs1}
\end{equation}
this implies that the flow preserves conformal class.  Next, we recall the transformation law
\begin{equation}
\Delta_{0} \phi - \frac12 R_{0} + \frac12 R e^{2\phi} = 0
\label{eq:trsc}
\end{equation}
relating the scalar curvatures $R_0$ and $R$ of two conformally related metrics, $g = e^{2\phi}g_0$. 
Using this, and writing $g(t) = u(t)g_0$, we see that \eqref{eq:rfs1} is equivalent to the scalar equation
\begin{equation}
\del_t u = \Delta_{g_0} \log u - R_0, \qquad u(0) = u_0,
\label{eq:rfs}
\end{equation}
and hence we focus on analysis on (\ref{eq:rfs}).

Combining \eqref{eq:trsc} and \eqref{eq:rfs1}, and noting that $\phi = \frac12 \log u$, we obtain the simple and
useful relationship
\begin{equation}
u_t = -Ru.
\label{eq:derivu}
\end{equation}

Now assume that $g(t)$ is a solution of \eqref{eq:rfs1} on some interval $0 \leq t < T$.  Differentiating \eqref{eq:trsc} yields
the evolution equation for scalar curvature,
\begin{equation}
\label{eq-R1}
\frac{\partial}{\partial t} R = \Delta R + R^2.
\end{equation}
It is also straightforward to check that
\begin{equation}
\frac{\partial}{\partial t} dA_t = -R dA_t,
\label{eq:evolaf}
\end{equation}
where $dA_t$ is the area form for $g(t)$.  

The Gauss-Bonnet formula can be extended to asymptotically conical surfaces by integrating $R$ over an increasing sequence of 
compact subdomains in $M$ with boundary a union of circles $\{r = \mbox{const.}\}$ on each end and then taking the limit; this yields
\begin{equation}
\int_M R\, dA_g = 4\pi\chi(M) - 4\pi \sum_{j=1}^\ell \alpha_j,
\label{eq:GB}
\end{equation}
where $\ell$ is the number of ends. In particular, it is important for our later considerations to note that 
if $\chi(M) \leq 0$ and if $M$ is AC, so that $\alpha_j > 0$, 
then $\int_M R < 0$.  For surfaces with finite total curvature which are not necessarily AC, one can define the  `aperture' of each end, and
there is an analogue of \eqref{eq:GB} due to Shiohama, see \cite{Wu}. 

The topological quantity $\chi(M)$ is obviously time independent. We can show that $\int_M R\, dV_g$ is constant along the flow as 
well, provided $R(\cdot,t) \in L^1(M, g(t))$ for each $t \geq 0$. Indeed, if we integrate (\ref{eq-R1}) on an expanding sequence
of sets as above, and if we use that $\del_r R = \calO(r^{-3-\tau})$,  we obtain
\begin{equation}
\label{eq-pres-R}
\frac{d}{dt}\int_M R\, dA  = \lim_{r \to \infty} \left(\int_{\partial B_0(p,r)} \del_\nu R\, d\sigma + \int_{B_0(p,r)} (R^2 - R^2)\, dA\right) = 0.
\end{equation}
This implies that $\sum_{j=1}^\ell \alpha_j$ is independent of $t$, but leaves open the possibility that for surfaces with more than one 
end, the individual cone angles might vary with $t$, so long as their sum remains constant.

\section{Existence}
The starting point for our discussion of \eqref{eq:rfs} is the following basic local existence result.
\begin{prop}
Let $(M,g_0)$ be a complete Riemannian surface with AC ends of order $\tau$. Then, for some $T > 0$, there is a unique solution 
to \eqref{eq:rfs} defined on $[0,T) \times M$, and each metric $g(t)$ in this family is AC of order $\tau$. 
\label{pr:ste}
\end{prop}

The existence of the flow solution $g(t)$ follows from the results and/or techniques appearing in several different papers; see in particular Wu \cite{Wu}, 
Schn\"urer, Schulze and Simon \cite{SSS} and Shi \cite{Sh}. Slightly less obvious is the fact that each of the metrics 
$g(t)$ remains asymptotically conical, and this is not discussed explicitly in any of these sources, although it does follow from an argument by 
Dai and Ma \cite{DM} which we sketch below. However, it is not difficult to give a direct proof of short-time existence for 
solutions of \eqref{eq:rfs} in the class of AC metrics by proving estimates for the corresponding linearized heat equation 
on appropriate weighted H\"older spaces and then constructing solutions of the nonlinear equation by a standard contraction 
mapping argument. We omit details, but we refer the reader to the analogous discussion in \cite{JMS}. 

We next turn to long-time existence for the solution. We wish to show
\begin{prop}
\label{prop-LTE}
Let $(M^2,g_0)$ be AC of order $\tau$. Then the solution to \eqref{eq:rfs} exists for all $t \geq 0$. Moreover, $\sup_{x \in M}|R(x, t)| 
\leq C$ for some fixed constant $C$ and for all $t \geq 0$.
\end{prop}

We follow the strategy set forth originally by Hamilton \cite{Ha} which employs a potential function $f$ for a metric $g$ to 
obtain uniform bounds on $R$, which leads to a priori estimates and the continuation of $u$ to all $t \geq 0$. 
By definition, such a function $f$ is a solution to the equation $\Delta_g f = R$ with bounded gradient. 
\begin{lem}
Let $(M,g)$ be an AC surface, and $R$ the scalar curvature of $g$. Then there exists a function $f$ on $M$ which satisfies
\[
\Delta_g f = R, \qquad \mbox{and} \qquad \sup_M |\nabla f| \le C. 
\]
\label{lem-potential}
\end{lem}

\begin{proof}
For an AC metric of order $\tau$, the scalar curvature $R_g$ lies in $r^{-2-\tau}\Lambda^{0,\alpha}$. 

We invoke some well-known Fredholm properties of the Laplacian on asymptotically conical spaces.
Define the $2\ell$-dimensional space (where $\ell$ is the number of ends)
\[
\mathcal{E} = \{\sum_{j=1}^\ell \chi_j(\beta_j \log r + \gamma_j): \beta_j, \gamma_j \in \RR\};
\] 
here each $\chi_j$ is a smooth function which equals one on the end $E_j$ and vanishes away from that end.
Then it is known that for each $0 < \nu < 1$, if $\nu \leq \tau$, the map 
\begin{equation}
\Delta_g: \mathcal{E} \oplus r^{-\nu}\Lambda^{2,\alpha}(M) \longrightarrow r^{-\nu-2}\Lambda^{0,\alpha}(M)
\label{eq:Fred}
\end{equation}
is surjective and its nullspace $\calN$ is $\ell$-dimensional.  (Implicit here is the easily verified fact that 
$\Delta_g (\log r) \in r^{-\nu-2}\Lambda^{0,\alpha}$.) This assertion is verified by combining Theorems 
4.20, 4.26 and 7.14 in \cite{M}. This route to the proof requires a lot of machinery, however, and especially 
in this low-dimensional setting, one can prove the same things in a more elementary way using barriers and 
sequences of solutions on compact exhaustions of $M$.

In any case, using \eqref{eq:Fred} with $\nu = \tau$, we obtain a solution $f \in \calE \oplus r^{-\tau}\Lambda^{2,\alpha}(M)$ 
to the equation $\Delta f = R$ which decomposes as
\[
f = \sum_{j=1}^\ell \chi_j (\beta_j \log r + \gamma_j) + \eta, \qquad \eta \in r^{-\tau}\Lambda^{2,\alpha}(M).
\]
for $\beta_j, \gamma_j \in \RR$. It is obvious from this that $\sup |\nabla f| \leq C$, so we are done.
\end{proof}

Every Riemann surface is a K\"ahler manifold, so instead of writing $g(t) = u(t)g_0$ we can alternately write $g_{i\bar{\jmath}} 
= (g_0)_{i\bar{\jmath}} + \partial_i\partial_{\bar{\jmath}}\phi$, where $\phi(t,\cdot)$ is called the K\"ahler potential of $g(t)$ 
relative to $g_0$. A standard calculation shows that if $g(t)$ satisfies the Ricci flow equation, then $\phi$ satisfies the Monge-Ampere equation
\begin{eqnarray}
\label{eq-pot}
\del_t \phi = \log\frac{\det((g_0)_{k\bar{\ell}} + \phi_{k\bar{\ell}})}{\det(g_0)_{k\bar{\ell}}} - 2f_0, \qquad \phi(x,0) =  0,
\end{eqnarray}
where $f_0$ is the potential function for the metric $g_0$ in the sense of Lemma \ref{lem-potential}, and $\phi_{k\bar{\ell}}$ 
is shorthand for $\partial_k \partial_{\bar{\ell}} \phi$.

Now define $f := -\frac{1}{2}\del_t \phi$. Differentiating \eqref{eq-pot} with respect to $t$ yields that $f$ satisfies the
evolution equation
\[
\del_t f = \Delta f, \qquad f(z,0) = -\frac12 \left. \del_t \phi \right|_{t=0} = f_0. 
\]
One can then show

\begin{claim}
The functions $f(z,t)$ defined above satisfy $\Delta f = R$  for every $t \ge 0$.
\end{claim}

\begin{proof}
To verify this claim, we carry through the following calculations.

First, we recall that 
\[
R g_{i\bar{\jmath}} = R_{i\bar{\jmath}} = -\partial_i\partial_{\bar{\jmath}}\log \det g_{i\bar{\jmath}}.
\]
Applying  $\partial_i\partial_{\bar{\jmath}}$ to (\ref{eq-pot}) yields
\[
-2\partial_i\partial_{\bar{\jmath}} f = - Rg_{i\bar{\jmath}} + R_0 (g_0)_{i\bar{\jmath}} - 2\partial_i\partial_{\bar{\jmath}} f_0;
\]
we then obtain by taking the trace with respect to $g_{i\bar{\jmath}}$ that
\[
\Delta f = R - \frac{1}{u}(R_0 - \Delta_0 f_0) = R;
\]
this uses both $\Delta_0 f_0 = R_0$ and $g = u g_0$. 

To conclude that $f$ is a potential function for $g(t)$, we still need to verify the boundedness of its gradient. 
\begin{lem}
If the solution $g(t)$ exists for $0 \leq t < T$, then $\sup_M |\nabla f(\cdot,t)|$ is finite for each $t < T$. 
\label{le:nablaf}
\end{lem}
\begin{proof}
We compute the evolution equation for $|\nabla f|^2$: since
\[
\frac{\partial\, }{\partial t}g^{ij} = -g^{ip}g^{jq}\frac{\partial\, }{\partial t}g_{pq} = R g^{ij},
\]
we have
\begin{equation}
\label{eq-nabla-f}
\begin{split}
\frac{\partial\, }{\partial t}|\nabla f|^2 &= (\frac{\partial\, }{\partial t}g^{ij})\nabla_if\nabla_j f + 2g^{ij}
\nabla_i\frac{\partial f}{\partial t}\nabla_j f \\
&= R |\nabla f|^2 + 2g^{ij}\nabla_i(\Delta f)\nabla_j f \\
&= R |\nabla f|^2 +2\langle\nabla R, \nabla f\rangle \\
&\le C_1|\nabla f|^2 + C_2t^{-1/2}|\nabla f|,
\end{split}
\end{equation}
where the last inequality uses Shi's estimates $|R(x,t)| \le C(t_1)$ and $|\nabla R| \le C(t_1)t^{-1/2}$ if the
flow exists on $M\times (0,t_1]$. 

Next, for each $x\in M$, fix $\tau = \tau(x)$ sufficiently small, and set $D(x,\tau) = \sup_{[0,\tau]} |\nabla f(x,t)|^2$. 
We know that $\sup_M D(x,0) < \infty$, so integrating (\ref{eq-nabla-f}) from $0$ to $\tau$ gives
\[
\begin{split}
D(x,\tau) & \leq C_0 + C_1 \tau D(x,\tau) + C_2 D(x,\tau) \int_0^{\tau} t^{-1/2}\, dt \\
& \leq C\left(1 + \tau D(x,\tau) + \sqrt{\tau D(x,\tau)}\right) \leq C(1+\sqrt{\tau} D(x,\tau)),
\end{split}
\]
assuming that $\tau < 1$ and $D(x,\tau) > 1$. The constant $C$ is independent of $x$, hence so is
$\tau = (4C^2)^{-1}$, and with this $\tau$ we obtain the uniform upper bound
\begin{equation}
\label{eq-first}
\sup_{M \times [0,\tau]} |\nabla f(x,t)|^2 \leq C.
\end{equation}

Finally, for $\tau \leq t < T$, 
\[
\frac{\partial\, }{\partial t}|\nabla f|^2 \le C_1|\nabla f|^2 +C_2 \tau^{-1/2}|\nabla f|
\]
so integrating from $\tau$ to any other value $t < t_0$ and using (\ref{eq-first}) gives
\[
\sup_M |\nabla f(\cdot,t)|^2 \le C(t)
\]
for any $t < T$, which is the desired result.
\end{proof}

We now complete the proof of Proposition \ref{prop-LTE}. Following \cite{Ha}, we define $h := \Delta f + |\nabla f|^2$ and the symmetric 
$2$-tensor $Z := \nabla^2 f - \frac{1}{2}\Delta f\cdot g$. A straightforward computation shows that 
\begin{equation}
\frac{\partial h}{\partial t} = \Delta h - 2|Z|^2.
\label{eq:evh}
\end{equation}
Using both Lemma \ref{le:nablaf} and $\Delta f = R$, we see that $h(\cdot, t)$ is bounded for each $t$. 
We can thus apply the maximum principle to (\ref{eq:evh}) to get that 
\[
\sup_M h(\cdot,t ) \le \sup_M h(x,0) \le C.
\]
This implies in turn that
\[
\sup_M R(\cdot,t) \le C,
\]
for all $t\in [0,T)$, where $C$ is independent of $T$.  In other words, the curvature remains uniformly bounded
for as long as the flow continues to exist.  Finally, since $u_t = -R u$, or equivalently, $(\log u)_t = -R$, 
we see that for some constants $C_1, C_2$ independent of $T$, $0 < C_1 \leq u \leq C_2$ for $0 \leq t < T$, 
and hence standard bootstrapping arguments show that $u$ remains bounded in $\calC^\infty$ for all $t < T$.
This continues to flow to $[0,T]$, and by Proposition \ref{pr:ste}, the flow continues to a slightly larger open interval.
This proves that $g(t)$ exists for all $0 \leq t < \infty$.
\end{proof}
 
\section{A priori bounds, I}
We now begin to examine the long-time behaviour of $g(t)$. In this section we obtain a number of estimates
concerning this behavior which we can prove using  variants of the maximum principle. These are not enough to complete the
proof of convergence, so in the next section we prove a number of further estimates using
quite different geometric arguments.  In all that follows, we often write $\Delta_t$ and $\nabla_t$ for
the Laplacian and gradient with respect to $g(t)$. 

The first result we obtain is that the asymptotically conical geometry is preserved. To prove this we closely follow an argument 
from Dai and Ma \cite{DM}.  
\begin{prop}
Suppose that $u_0(x) \in r^{-\tau}\Lambda^{2,\alpha}(M)$ and let $u(x,t)$ be the solution to \eqref{eq:rfs} with
$u(x,0) = u_0(x)$. Then for all $t \geq 0$, $u(\cdot,t) \in r^{-\tau}\Lambda^{2,\alpha}(M)$, with bounds uniform
in any strip $0 \leq t \leq T$.
\label{pr:acpersist}
\end{prop}
\begin{proof}
Let $r$ denote a smooth function which agrees with the radial distance function of the model conic metric on each end, and 
such that $r \geq 2$ on all of $M$. We may choose  $r$ so that it satisfies $|\Delta_t r| + |\nabla_t r| \leq C(t)$.

We first show that $|R(t)| \leq C(t) r^{-2-\tau}$. Using \eqref{eq-R1} and the absolute bound $|R(t)| \leq C$, we derive that
\[
\del_t (R^2) \leq \Delta_t R^2 - 2|\nabla_t R|^2 + 2C R^2,
\]
from which we calculate that $w = r^{4 + 2\tau}R^2$ satisfies
\[
\del_t w \leq \Delta_t w + A \cdot \nabla_t w + B w,
\]
where $A$ and $B$ are uniformly bounded. Hence, for $C_1, C_2$ sufficiently large, 
\[
\del_t (w - C_1 e^{C_2 t}) \leq \Delta_t (w - C_1 e^{C_2 t}) + A \cdot \nabla_t (w - C_1 e^{C_2 t} ),
\]
and moreover, $w - C_1 e^{C_2 t} \leq 0$ at $t=0$.  We now invoke the maximum principle for an
evolving family of metrics, proved by Ecker and Huisken and recorded as Theorem \ref{thm-EH} in the appendix here.
This gives 
\[
|R| \leq C_1' e^{C_2' t} r^{-2-\tau}.
\]
Integrating $ (\log u)_t = -R$ from $0$ to $t$ yields
\[
\log u(x,t) - \log u(x,0) = \int_0^t R(x,s)\, ds,
\]
which completes the proof. 
\end{proof}

\begin{cor}
Let $(M^2,g_0)$ be AC of order $\tau$. Then $(M,g(t))$ remains AC of order $\tau$ for each $t \geq 0$. 
\end{cor}

The next two results concern upper and lower bounds for $u$ which are global in $t$.

The linear upper bound for $u(t)$ follows from a well-known argument due to Aronson and Benilan \cite{AB}.
\begin{prop}
\label{pr:ab}
Let $u$ be a solution of \eqref{eq:rfs}. Then $u_t \le \frac{u}{t}$, and hence
\[
u(t) \leq C(1+t)
\]
for some constant $C>0$ and all $t \geq 0$. 
\end{prop}
\begin{proof} Following \cite{AB}, we define $u_{\lambda}(x,t) := \lambda u(x, \lambda^{-1}t)$. This satisfies 
\[
\frac{\del u_\lambda }{\del t} = \Delta_{g_0} \log u_{\lambda}  - R(g_0), \qquad  \mbox{and} \qquad
\left. u_{\lambda}(x,t)\right|_{\lambda = 1} = u(x,t).
\]
Moreover, for $\lambda > 1$, $u_{\lambda}(x,0) = \lambda u(x,0) > u(x,0)$. Setting $v_{\lambda}(x,t) := u_{\lambda}(x,t) - u(x,t)$, 
we calculate that
\begin{equation}
\label{eq-dif}
\frac{\del v_\lambda}{\partial t} (x,t) = \Delta_{g_0} (a(x,t)\cdot v_{\lambda}), \qquad v_{\lambda}(x,0) \geq 0, 
\end{equation}
where 
\[
a(x,t) := \int_0^{1}\frac{d\theta}{\theta u_{\lambda}(x,t) + (1-\theta)u(x,t)}\, d\theta.
\] 

Using Proposition \ref{pr:acpersist}, for any $T > 0$ and all $(x,t) \in M \times [0,T]$, $C_1(T) \leq u(x,t) \leq C_2(T)$; hence $a(x,t)$ 
is also bounded above and below by ($T$-dependent) constants. This shows that \eqref{eq-dif} is strictly parabolic on any finite time 
interval. Since $v_{\lambda}$ is uniformly bounded on $M \times [0,T)$ and since $v_{\lambda}(x,0) \ge 0$, it follows from  the maximum principle that $v_\lambda \geq 0$ for $t \geq 0$.  We also have $v_1(x,t) = 0$, so we find that   $\lambda \mapsto v_\lambda$ is nondecreasing  in some interval $[1,\lambda_0)$; that is, $\left. \del_\lambda \right|_{\lambda = 1} v_\lambda(x,t) 
\geq 0$, or equivalently, $u_t \le t^{-1} u$. Finally, $u > 0$, so the other estimate in the statement of this result follows by
integrating $u_t/u \leq 1/t$. 
\end{proof} 

It is also not hard to prove a uniform lower bound for $u$. 
\begin{prop}
There exists a constant $C_1 > 0$ such that $C_1 \leq u(x,t)$ for all $(x,t) \in M \times [0,\infty)$. 
\label{pr:clb}
\end{prop}
\begin{proof}
The first step is to show that the initial metric is conformal to another AC metric (with the same cone angles on each end) 
with $R_0 \leq 0$. Not surprisingly, this relies on the assumption that $\chi(M) < 0$.

Recall that if $\hat{g}_0 = e^{2\psi} g_0$, then $\Delta_0 \phi - \frac12 R_0 + \frac12 \hat{R}_0 e^{2\phi} = 0$ (where
we denote the scalar curvatures of $g_0$ and $\hat{g}_0$ by $R_0$ and $\hat{R}_0$, respectively). Thus, given an AC metric $g_0$, to show that there exists $\psi$ such that $\hat{g}_0$ is AC and has scalar curvature $\hat{R}_0 \le 0$, it is sufficient to obtain $\psi$ satisfying  $\Delta_0 \psi \geq \frac12 R_0$. 

Let us choose a function $Q \in r^{-2-\tau}\Lambda^{0,\alpha}(M) \cap \calC^\infty$ which satisfies $Q(x) \geq \frac12 R_0(x)$ and 
$\int Q = 0$. This is possible since, by \eqref{eq:GB}, $\int R_0 < 0$.  As a consequence of the surjectivity of the map $\Delta_g$ in \eqref{eq:Fred}, we can find $\psi_1 \in \calE \oplus 
r^{-\tau}\Lambda^{2,\alpha}$ with $\Delta_0 \psi_1 = Q$. This solution may grow logarithmically, so we must modify it further by 
adding on a function $w$ in the nullspace of the mapping \eqref{eq:Fred} such that the coefficient $\beta_j$ of $\log r$ in the 
expansion of $w$ on each end $E_j$ equals the corresponding coefficient of $\log r$ in the expansion for $\psi_1$ on that end. 
Proposition 6 in \cite{JMS} shows that this is possible.  Therefore $\psi = \psi_1 - w$ is bounded and satisfies all the required properties. 

Now we write the evolving metric  in the form $g(t) = u(t) g_0 = u(t) e^{-2\psi}\hat{g}_0 := \hat{u}(t) \hat{g}_0$. Then $\hat{u}(0) = u(0)\cdot e^{-2\psi} \geq C_1 > 0$ and, using
that $\hat{R}_0 \leq 0$, 
\[
\del_t \hat{u} = \Delta_{\hat{g}(0)} \log \hat{u} - \hat{R}_0 \geq \Delta_{\hat{g}(0)} \log \hat{u}.
\]
Hence by the minimum principle, $\hat{u}(t) \geq \hat{u}(0) \geq C_1 > 0$ for all $t \geq 0$. The lower bound for $u(x,t)$ now follows from that for $\hat{u}$, together with the relation $u=\hat{u} e^{2\psi}$.
\end{proof}

Since it is also proved using the maximum principle, we include one final result, that if the initial curvature is nonpositive,
then the curvature remains nonpositive for all time. Interestingly, our main convergence result is significantly easier to prove under
this assumption on $R_0$. 

\begin{prop}
Let $R_0 \leq 0$. Then $R \leq 0$ for all $t \geq 0$. 
\end{prop}
\begin{proof}
Recall that $\del_t R = \Delta_0 R + R^2$. For any $A > 0$ define 
\[
Q(x,t) = R(x,t) + \frac{1}{A+t}.
\]
Then
\[
\del_t Q = \del_t R - \frac{1}{(A+t)^2} = \Delta_0 Q + R^2 - \frac{1}{(A+t)^2} = 
\Delta_0 Q + VQ,
\]
where
\[
V = R - \frac{1}{A+t}.
\]
Since $R(x,0) \leq 0$, the function $V$ is strictly negative for $t=0$.  For fixed $A$, we now define
\begin{multline*}
T = T_A = \inf \{ t > 0: R(x,t) < \frac{1}{A+\tau}\ \mbox{for}\ 0 < \tau < t \\ 
 \mbox{and} \ R(x_0,t) = \frac{1}{A+t} \ \mbox{for some}\ x_0 \in M\}.
\end{multline*}
Note that, since $R(x,0) \leq 0$, $T$ is strictly positive, and furthermore, if $T= \infty$, then
$R(x,t) \leq 1/(t+A)$ for all $t \geq 0$.

An application of the maximum principle to $Q(x,t)$ on $M \times [0,T)$, using that $V \leq 0$ on this domain, gives
\[
Q(x,t) \leq \min \{ \inf_{t \in [0,T]} \frac{1}{t+A}, \frac{1}{A} \} \leq \frac{1}{A}
\]
for all $(x,t)$ in this domain, or in other words,
\[
R(x,t) \leq \frac{1}{A} - \frac{1}{t+A} = \frac{t}{A(t+A)}.
\]
Since $t/(t+A)$ is increasing in $t$, we have
\[
R(x,t) \leq \frac{T}{A(T+A)} = \frac{1}{T+A} \frac{T}{A}.
\]
However, we know that $R(x_0,T) = 1/(T+A)$ for some $x_0$, which implies that $T_A \geq A$.

These calculations show that if $A > 0$ is arbitrary, then either $R(x,t) \leq 1/(t+A)$ for all $t \geq 0$,
or else $R(x,t) \leq \frac{t}{t(t+A)}$ for $0 \leq t \leq A$.  Letting $A \to \infty$ in either case implies
that $R(x,t) \leq 0$ for all $t \geq 0$. 
\end{proof}

\section{A priori bounds, II}
Since $g(t)$ has a uniform linear upper bound, it is natural to consider the rescaled family of metrics

\[
\tilde{g}(t) := \frac{1}{t} g(t);
\]
indeed, it follows from Proposition \ref{pr:ab} that  $\tilde{g}(t) \leq C g(0)$ for all $t \geq 1$, with $C$ independent of $t$. This 
family also satisfies an evolution equation: setting $\tau = \log (t)$, we calculate that
\begin{equation}
\del_\tau \tilde{g}(\tau) = - (\tilde{R} + 1) \tilde{g}(\tau),
\label{eq:evgtilde}
\end{equation}
where $\tilde{R}$ is the scalar curvature of $\tilde{g}$ at time $\tau$, or equivalently, 
\begin{equation}
\del_\tau \tilde{u} = \Delta_0 \log\tilde{u} - \tilde{u} - R_0 = - (\tilde{R} + 1)\tilde{u},
\label{eq:evutilde}
\end{equation}
where $\tilde{u}(\cdot,\tau) = u(\cdot, t)/(t) = e^{-\tau} u(\cdot, e^\tau)$.  
\begin{prop}
The function $\tilde{u}(\tau,x)$ is monotone nonincreasing in $\tau$ for each fixed $x$.
\end{prop}
\begin{proof} One form of the original evolution equation is that $u_t = -Ru$; on the other hand,
by Proposition \ref{pr:ab}, 
\[
R(x,t) \ge -\frac{1}{t} \Rightarrow  \tilde{R}(\tau) \geq -1.
\]
From these we see that the right side of (\ref{eq:evutilde}) is nonpositive;  hence $\del_\tau \tilde{u} \leq 0$.
\end{proof}

We now state the main result of this paper. 
\begin{thm}
\label{thm-part}
The metric $\tilde{g}(\tau)$ converges in $\calC^\infty$ on every compact set as $\tau\to\infty$; the limiting metric $g_\infty$ 
is complete,  hyperbolic, and has finite area. 
\end{thm}

The proof will occupy the rest of this section. It proceeds roughly as follows: since $\tilde{u}$ is monotone nonincreasing in $\tau$, 
it has a limit as $\tau \nearrow \infty$. If we can show that $\tilde{u}(x_0,\tau)\geq c > 0$ for some fixed $x_0$, then a gradient estimate for $\tilde{u}$ together with a Harnack inequality for $R(\tilde{g})$
imply that $\tilde{u}$ stays bounded away from zero in a fixed neighbourhood of $x_0$. This allows us to apply Hamilton's compactness 
theorem to solutions of this flow and thereby conclude that $(M,\tilde{g}(\tau), x_0)$ converges to a {\it complete} Riemannian surface. 

The first step, that $\tilde{u}(x,\tau)$ cannot tend to $0$ for every $x$, is accomplished using a further rescaling of the metric $g$, and 
requires the topological hypothesis that $\chi(M) < 0$. 
\begin{prop}
There exists a point $x_0 \in M $ such that $\tilde{u}(x_0,\tau) \geq \delta > 0$ for some fixed $\delta$ and all $\tau \ge 0$.
\label{pr:lowerboundtu}
\end{prop}
\begin{proof}
By monotonicity of $\tilde{u}$ in $\tau$, $\lim_{\tau \to \infty} \tilde{u}(x,\tau) := \tilde{U}(x)$ exists for every $x \in M$. 
Thus we must prove that $\tilde{U} \not\equiv 0$. 

Suppose, to the contrary, that $\tilde{U}(x) = 0$ for every $x$. Dini's theorem states that a monotone sequence of
continuous functions which converges pointwise to a continuous function must in fact converge uniformly on
compact sets. We show that this leads to a contradiction.

Let $K$ be a sufficiently large compact set so that $M \setminus K$ is a union of (asymptotically conical) ends 
$E_j \cong [1,\infty) \times S^1$ (using the coordinates $(r,\theta)$ where the metric has the form
$dr^2 + \alpha^2 r^2 d\theta^2 + \calO(r^{-\nu})$).  Choose any sequence $\tau_i \nearrow \infty$ and $p_i \in K$ so that 
\[
\alpha_i := \tilde{u}(p_i,\tau_i) = \max_{x\in K} \tilde{u}(x,\tau_i).
\]
By hypothesis, $\alpha_i \to 0$.

Now perform yet another rescaling: set
\[
\bar{g}_i(\tau) = \alpha_i^{-1} \tilde{g}(\tau_i + \alpha_i\tau),
\]
and let $\bar{u}_i$ be the corresponding conformal factor. Note that $R(\bar{g}_i) = \alpha_i R(\tilde{g}_i)$, and by construction, 
$\bar{u}_i(p_i, 0) = 1$. 

Since $R(\tilde{g}_i) \geq -1$, we have a lower bound $R(\bar{g}_i) \geq -\alpha_i \geq -C$, but we do not yet know that
the curvatures of this sequence of metrics are uniformly bounded from above. We first prove our result assuming that such
an upper bound is true; i.e.\ that
\begin{equation}
R(\bar{g}_i) \leq C. 
\label{eq:Rbar0}
\end{equation}
Then at the end we justify (\ref{eq:Rbar0}), allowing the constant $C$ in (\ref{eq:Rbar0}) to depend on a compact set over which we are estimating the rescaled curvature; this  will be enough to finish the argument.  

The immediate goal is to estimate $\nabla \bar{u}_i$, which will allow us to prove that $\bar{g}_i$ converges to
a complete metric.

Let $f_0$ be the potential function associated to the original scalar curvature function $R_0$, as defined in Lemma \ref{lem-potential},
and let $f(x,t)$ be its evolution under the linear heat flow $\del_t f = \Delta_t f$.  It has been proved   in \S 3 that $f(x,t)$ is a potential function for $g(t)$.  We claim that $\log u(x,t) \equiv f_0(x) - f(x,t)$, or equivalently, that $k(x,t) := \log u(x,t) - f_0(x) 
+ f(x,t) \equiv 0$.  To prove this, note that $\del_t k = \Delta_t k$ and $k(x,0) = 0$. Furthermore, $x \mapsto k(x,t)$ is
bounded for each $t$, which holds because $|\del_t f| = |\Delta f| = |R|$ is bounded for each $t$. Hence $|f(x,t) - f(x,0)| \leq Ct$,
and $|\log u|$ is bounded for each $t$ as well. The maximum principle then implies that $k \equiv 0$. 

Continuing on, we have
\[
\frac{\partial}{\partial\tau}(\log \tilde{u} - f_0) = \Delta_{\tilde{g}}(\log\tilde{u} - f_0) - 1,
\]
and then  a computation from \cite{Ha} gives us
\[
\frac{\partial}{\partial\tau}|\nabla(\log\tilde{u} - f_0)|_{\tilde{g}}^2 \le \Delta_{\tilde{g}}|\nabla(\log\tilde{u} - f_0)|_{\tilde{g}}^2.
\]
Since $|\nabla(\log\tilde{u} - f_0)|_{\tilde{g}}^2(\cdot,\tau)$ is bounded for each $\tau$, the maximum principle shows that
\begin{equation}
\label{eq-nabla-100}
|\nabla_{\tilde{g}} \log \tilde{u}|_{\tilde{g}}(\cdot,\tau) \le C \qquad \mbox{for all} \quad \tau \geq 0,
\end{equation}
or equivalently,
\[
|\nabla_{\bar{g}_i} \log\bar{u}_i|_{\bar{g}_i} \le C\sqrt{\alpha_i},
\]
and therefore (since $\bar{u}_i(x,0) \le 1$ on $K$)
\[
|\nabla\bar{u}_i(\cdot, 0)|_{\bar{g}_i} \le C\sqrt{\alpha_i} \qquad \mbox{on} \quad K.
\]
Combining this with the fact that $\bar{u}_i(p_i,0) = 1$, we see that there exists an $\eta > 0$ such that for each sufficiently
large $i$, 
\[
\bar{u}_i(x,0) \geq \eta \qquad \mbox{for all} \quad x\in B_{\bar{g}_i(0)}(p_i,1);
\]
from this we also see that
\begin{equation}
\label{eq-inj-bari}
\inj_{\bar{g}_i(0)}(p_i) \ge \delta > 0 \qquad \mbox{for all} \quad i.
\end{equation}

Boundedness of $R(\bar{g}_i)$  and the estimate (\ref{eq-inj-bari}) are precisely the hypotheses needed to apply Hamilton's 
compactness theorem. This result states that, after passing to a subsequence, $(M,\bar{g}_i(\tau),p_i)$ converges 
in the pointed Cheeger-Gromov sense to a limiting family of complete metrics $(M_{\infty}, g_{\infty}(\tau), p_{\infty})$ which is an 
eternal solution of the Ricci flow. 
Now, it follows from the Aronson-B\'enilan inequality that any ancient solution of the Ricci flow on surfaces is nonnegatively 
curved. Moreover, the work of B.L.\ Chen \cite{Chen} shows that the scalar curvature of any ancient solution is nonnegative. 
However, since $\tilde{u} (\tau_i) \leq \alpha_i$ on $K$, we have that
\[
\mbox{diam}_{\tilde{g}(\tau_i)} (K) \leq C \sqrt{\alpha_i} \Longrightarrow \mbox{diam}_{\bar{g}_i}(K) \leq C',
\]
where the last constant is independent of $i$. Hence $\bar{g}_i$ is $\calC^2$ close to $g_\infty$ on the compact set $K$. But this leads to a 
 contradiction since it follows from (\ref{eq:GB}) that $\limsup_{i \to \infty} \int_K R(\bar{g}_i) < 0$. 

It remains to verify (\ref{eq:Rbar0}).  

Let $\beta_i = \alpha_i \tilde{R}(p_i, \tau_i + \alpha_i)$ (this is the curvature of $R(\bar{g}_i)$ at $p_i$ at $\tau = 1$). We claim that
$\beta_i$ is bounded.  If this were to fail, i.e.\ if $\beta_i \to \infty$, at least along some subsequence, then by the Harnack 
inequality for $\tilde{R}$ (the form of this inequality which we use here is stated in \cite{Ch}; the original Harnack estimate for Ricci flow, proven by Hamilton, appears in \cite{Ha2}), 
 for $\alpha_i \leq \tau \leq \alpha_i + 1$, we would have
\[
\tilde{R}(p_i,\tau_i+\tau) \geq C (\tilde{R}(p_i,\tau_i+\alpha_i) + 1)\cdot e^{-C(\tau - \alpha_i)} - 1 \geq C' \frac{\beta_i}{\alpha_i},
\]
with $C'$ independent of $i$. Moreover,
\[
\del_\tau\tilde{u} = -(\tilde{R} + 1)\tilde{u} \leq - C' \beta_i
\]
in this interval, so that
\[
\tilde{u}(p_i,\tau_i+\alpha_i + 1) \leq \tilde{u}(p_i,\tau_i+\alpha_i) - C' \beta_i.
\]
Using monotonicity and iteration, we have $\tilde{u}(p_i,\tau_i + \alpha_i) \to -\infty$ as $i \to \infty$. This 
contradicts the fact that $\tilde{u} > 0$. 

We need to apply Hamilton's compactness theorem again as in the argument on the previous page; in fact, we need a generalization, proved in Appendix E of \cite{KL},
which requires that the curvature $R(\bar{g}_i)$ be bounded uniformly in time, but only over a fixed compact set $K \subset M$. 
This hypothesis is verified as follows.

The Harnack estimate for $\tilde{R}$ states that for $0 \leq s \leq 1/2$, 
\[
\tilde{R}(p,\tau_i+s\alpha_i) + 1 \le C(1+\tilde{R}(p_i,\tau_i+\alpha_i)) \, \exp \left(
C\frac{\dist^2_{\tilde{g}(\tau_i+s\alpha_i)}(p,p_i)}{(1-s)\alpha_i}\right),
\]
and hence
\[
\tilde{R}(p,\tau_i+s\alpha_i)  \leq \frac{C(\rho)}{\alpha_i},
\]
for all $p$ with $\dist_{\tilde{g}(\tau_i)}(p,p_i) \le \sqrt{\alpha_i} \, \rho$. Here we are use the fact that the distance with respect to 
$\tilde{g}(\tau)$ between any two fixed points decreases in $\tau$; this follows from the monotonicity of $\tilde{u}$. 
These facts together imply the bound we are seeking, that 
\[
|\bar{R}_i(s,x)| \le C(\rho) \qquad x \in B_{\bar{g}_i}(p_i,\rho),
\]
since $s \leq 1/2$.  

The other hypothesis we must verify is that $\bar{g}_i(0)$ has a bound on its injectivity radius at $p_i$ which is uniform in  $i$.
This is done exactly as before. 

The compactness theorem allows us to pass to a limiting metric, which is necessarily an ancient solution of the Ricci flow.  As 
we have already discussed, any ancient solution has nonnegative curvature, but considering the Gauss-Bonnet integral over $K$,
we obtain the same contradiction as before.

This completes the proof that $\alpha_i$ must remain bounded away from zero. 

\begin{claim}
There exists an $x_0\in K$ so that $\tilde{u}(x_0,\tau) \ge \delta > 0$ for all $\tau \ge 0$.
\end{claim}

\begin{proof}
Since all $p_i \in K$, the set $\{p_i\}$ has an accumulation point in $K$, call it $x_0$, with the property that
\begin{equation}
\label{eq-closeness}
\lim_{i\to\infty}\dist_{\tilde{g}(0)}(p_i,x_0)  = 0.
\end{equation}
Since $\alpha_i$ is bounded, there exists $\delta > 0$ so that
$$\tilde{u}(p_i,\tau_i) \ge 2\delta, \,\,\, \mbox{for all} \,\,\, i.$$ 
Combining (\ref{eq-nabla-100}) with the uniform global upper bound on $\tilde{u}$, we have $|\nabla\tilde{u}| \le C$. This implies that 
$$\tilde{u}(x_0,\tau_i) \ge \tilde{u}(p_i,\tau_i) - C\dist_{\tau_i}(p_i,x_0) \ge 2\delta - C\dist_{\tilde{g}(0)}(x_0,p_i) \ge \delta,$$
for $i$ sufficiently big so that $\dist_{\tilde{g}(\tau_i)}(x_0,p_i) \le \dist_{\tilde{g}(0)}(x_0,p_i) \le \frac{\delta}{C}$. 

(Here we have used  (\ref{eq-closeness}) together with the result that the distances are nonincreasing in $\tau$.)
\end{proof}

With this claim, we complete the proof of Proposition \ref{pr:lowerboundtu}.

\end{proof}

We now know that there exists some $x_0\in M$ such that $\tilde{u}(x_0,\tau) \geq c > 0$ for all $\tau\geq 0$.
The gradient estimate (\ref{eq-nabla-100}) can be applied as before to show that for some $r_1 > 0$,
\begin{equation}
\label{eq-inj}
\tilde{u}(x,\tau) \ge \delta > 0
\end{equation}
for all $x\in B_{\tilde{g}(\tau)}(x_0,r_1)$, $\tau\geq 0$.  Clearly $\tilde{g}(0)$ is `$\kappa$-noncollapsed'; i.e.,\ there exist constants $\kappa, r_0 > 0$ ($r_0 \le r_1$) so that if $B_0(x,r)$ is the geodesic ball 
around $x$ of radius $r \le r_0$, with respect to $\tilde{g}(0)$, then $\mathrm{Area}_{\, 0}\,(B_0(x,r)) \geq \kappa r^2$.
If $B_\tau(x,r)$, for $r \le r_0$, is the corresponding geodesic ball with respect to $\tilde{g}(\tau)$, then by monotonicity of $\tilde{u}$,
$B_0(x_0,r) \subset B_{\tau}(x_0,r)$. Hence using (\ref{eq-inj}), we have
\begin{equation}
\label{volume-lower}
\begin{array}{rcl}
\mathrm{Area}_{\, \tau}(B_{\tau}(x_0,r)) &=& \int_{B_{\tilde{g}(\tau)}(x_0,r)} \frac{\tilde{u}(x,\tau)}{\tilde{u}(x,0)}\, dV_{\tilde{g}(0)}  \\
&\ge& \tilde{\delta}\int_{B_{\tilde{g}(\tau)}(x_0,r)} dV_{\tilde{g}(0)} \ge 
\int_{B_0(x_0,r)} dV_{\tilde{g}(0)} \\
&=& \mathrm{Area}_{\, 0} (B_0(x_0,r))\ge \kappa' r^2.
\end{array}
\end{equation}
In other words, $\tilde{g}(\tau)$ is $\kappa'$-noncollapsed, where $\kappa'$ is a $\tau$-independent multiple of $\kappa$.

We shall use the same compactness theorem as before to show that $\tilde{g}(\tau)$ converges to a complete metric. 
In order to do so, we must show that $|R(\tilde{g}(x,\tau)| \leq C(\rho) < \infty$ for $x \in B(x_0,\rho)$, for all $\rho > 0$  and that  $\mathrm{Inj}_{\tilde{g}(\tau)}(x_0) \geq c > 0$.

For the next step, we establish a local curvature bound. 
\begin{lem}
\label{lem-cur-comp}
For every $\rho > 0$ there is a constant $C_\rho$ so that 
\[
\tilde{R}(x,\tau) \le C_\rho \qquad \mbox{for all} \quad (x,\tau) \ \mbox{such that}\ x \in B_{\tau}(x_0,\rho),\  \tau \ge \tau_0.
\]
\end{lem}
\begin{proof}
Using the Harnack inequality (\ref{eq-har-r}) quoted in the appendix, it suffices to show that for some $C > 0$, 
$\tilde{R}(x_0,\tau) \le C$ for all $\tau \ge 0$.  If this were to fail, then there would exist a sequence $\tau_i\nearrow \infty$ 
for which $Q_i := \tilde{R}(x_0,\tau_i) \to \infty$. Taking $\tau_2 = \tau_i + \tau$, $\tau_1 = \tau_i$,  $x_2 = x_1 = x_0$ 
in (\ref{eq-har-r}), and fixing $A > 0$, we see that for  $0 \leq \tau \leq A$, 
\begin{equation}
\label{eq-R-bound1}
\tilde{R}(x_0,\tau_i+\tau) \geq e^{-C\tau}\cdot (\tilde{R}(x_0,\tau_i) + 1) - 1  \geq cQ_i,
\end{equation}
where $c$ depends only on $A$. 

Arguing as before, $\del_\tau \log\tilde{u}(x_0,\tau) \le -c Q_i \to -\infty$, so integrating over the interval $[\tau_i,\tau_i+\tau]$
(with $\tau \le A$), we obtain
\[
\tilde{u}(x_0,\tau_i+\tau) \le e^{-c Q_i\tau}\tilde{u}(x_0,\tau_i) \to 0,
\]
which contradicts that $\tilde{u}(x_0,\tau) \geq c > 0$.

This proves that $\tilde{R}(x_0,\tau) \le C$ for all $\tau \geq 0$. Finally, putting $x_1 = x$, $x_2 = x_0$, $\tau_1 = \tau$ and 
$\tau_2 = \tau + 1$ in (\ref{eq-har-r}), we obtain the conclusion of the Lemma.
\end{proof}

Finally, we invoke a result due to Cheeger, who has shown that for a geodesic ball $B(x_0,r_0)$ in some Riemannian manifold $(M^n,g)$,
if there are lower bounds on Ricci curvature and volume, then the injectivity radius of $(M,g)$ at $x_0$ is bounded away from $0$
by some constant depending only on these bounds. We have established all of these hypotheses (in Lemma \ref{lem-cur-comp} and (\ref{volume-lower})); hence
\begin{equation}
\label{eq-injradest}
\inj_{\tau}(x_0) \geq c > 0 \qquad \forall \, \tau \geq 0.
\end{equation}

We now apply the compactness theorem to conclude that for every sequence $\tau_i\to \infty$, there is a subsequence of 
$(M,\tilde{g}(\tau_i + \tau), x_0)$ which converges smoothly as a family of pointed spaces to a smooth complete family of 
metrics $(M_\infty, g_{\infty}(\tau),x_\infty)$.  More specifically, for each 
compact interval $I\subset [0,\infty)$ and for any compact set $K\subset M_{\infty}$ containing $x_0$, 
there are pointed, $\tau$-independent diffeomorphisms $\phi_{K,i}:K\to K_i \subset M$ such that the appropriate 
subsequence of $\phi^*_{K,i}\tilde{g}(\tau_i+\tau)$ converges smoothly to $g_{\infty}(\tau)$ on $K\times I$. 

It is clear from the monotonicity of $\tilde{g}(\tau)$ that this limit is unique and does not depend on $\tau$, so $g_\infty$ satisfies 
the stationary equation.  Thus $\tilde{u}(\cdot,\tau)$ converges uniformly on compact sets of $M$ to a continuous function $\tilde{U}(x)$.
This limiting function is nonnegative, but there is still a possibility that it vanishes on some nontrivial closed set, 
which would make $\tilde{U}(x) g_0$ degenerate. However, this is ruled out by uniqueness of the limiting metric and the fact that 
$g_\infty$ is complete, so $\tilde{U} > 0$ everywhere. 

We may now use the equation
\[
\Delta\log \tilde{u} - R_0 = -\tilde{R}\tilde{u}.
\]
and the fact that $|\log \tilde{u}| \leq C$ and $-1 \le \tilde{R}(x,\tau) \le C(K)$ for $\tau \ge 0$ and $x\in K$ to obtain
$\tau$-independent bounds on all higher order derivatives of $\tilde{u}$ over any compact set by standard bootstrapping. It follows that the convergence of $\tilde{u}(x, \tau)$ to $\tilde{U}(x)$ is
 $\calC^\infty$ on compact sets. We have thus proved that 
\begin{equation}
\label{eq-smooth-conv}
\tilde{g}(x,\tau) = \tilde{u}(x,\tau) g_0 \longrightarrow \tilde{U}(x) g_0 = g_\infty
\end{equation}
smoothly on compact sets. 

It remains to show that $g_\infty$ is hyperbolic and has finite area. For the first of these, recall that we already know that 
$R_{\infty}(x) \ge -1$. If we assume that $R_{\infty}(y) \ge -1 + 2\delta$ for some $\delta > 0$ and some $y\in B_{g_{\infty}}(x_0,r_0)$, where
$x_0$ and $r_0$ are as above, then $\tilde{R}(y,\tau) \ge -1 +  \delta$ for $\tau \ge \tau_0$ and $y\in B_{\tilde{g}(\tau)}(x_0,r_0)$ and it follows that
\[
\del_\tau \log\tilde{u}(y,\tau) = -(\tilde{R} + 1) \le -\delta \Longrightarrow 
\tilde{u}(y,\tau) \le \tilde{u}(y,\tau_0) \cdot e^{-\delta(\tau - \tau_0)},
\]
which contradicts (\ref{eq-inj}). Therefore, $R_{\infty} \equiv -1$ on $B_{g_{\infty}}(x_0,r_0)$. Now put $\tau_1 = \tau$ and $\tau_2 = \tau + 1$ in 
(\ref{eq-har-r}) and use (\ref{eq-delta}) to get that
\begin{equation}
\label{eq-previous}
\tilde{R}(x_1,\tau) + 1 \le \tilde{R}(x_2,\tau+1) e^{d(x_1,x_2,\tau)^2/4 + C}.
\end{equation}
Taking $x_2\in B(x_0,r)$ and $x_1$ any other point in $M$, and letting $\tau\to\infty$ in (\ref{eq-previous}), the smooth convergence
of $\tilde{u}$ implies that
\[
\tilde{R}_{\infty}(x_1) + 1 \leq \tilde{R}(x_2) e^{d_{\infty}(x_1,x_2)^2/4 + C} = 0.
\]
However, since $\tilde{R}_{\infty} + 1 \ge 0$, we conclude finally that $R_{\infty} \equiv -1$ on $M_\infty$.

The final step concerns the finiteness of the area. Since $\tilde{g}_\infty$ is hyperbolic and complete, if its area were not finite,
then its area growth would be exponential; i.e.,\ there would exist constants $C_1, C_2, a, b > 0$ so that 
\[
C_1 e^{ar} \le \mathrm{Area}_\infty (B_{\infty}(p,r)) \le C_2 e^{br}
\]
for $r$ large. This would imply that for $i$ large, 
\begin{equation}
\label{eq-hyper-growth}
\frac12 C_1 e^{ar} \le \mathrm{Area}_{\tau_i} B_{\tau_i}(p_i,r) \le 2C_2 e^{br}.
\end{equation}
However, $\tilde{u}(\cdot, \tau_i) \leq c_2$ uniformly on $M$, for all $i$, so $(M,\tilde{g}(\tau_i))$ can have at most
quadratic area growth. This contradiction finishes the proof of Theorem \ref{thm-part}.

\section{Appendix: maximum principle and Harnack inequality}
\label{sec-maximum}

We state two versions of the maximum principle for complete manifolds, and then a version of the Harnack estimate which holds for such geometries.

\begin{thm}[Maximum principle]
\label{thm-max}
Let $g(t)$, $0 \le t < T$, be a family of complete Riemannian metrics on a noncompact manifold $M$ which
vary smoothly in $t$ and which satisfy $C_1 g(0) \leq g(t) \leq C_2 g(0)$ for some fixed constants $C_1, C_2$
and for all $t \in [0,T)$.  Let $f(x,t)$ be a smooth bounded function on $M\times [0,T)$ which satisfies the initial condition 
 $f(x,0) \geq 0$, and satisfies the parabolic equation 
\[
\frac{\del \,}{\del t}f = \Delta_{g(t)} f + Q(f,x,t),
\]
where $Q(f,x,t) \ge 0$ whenever $f \le 0$. Then $f(x,t) \ge 0$ on $M\times [0,T)$.
\end{thm}
The proof can be found in \cite{Sh} and \cite{Sh1}.


\begin{thm}[Ecker-Huisken Maximum principle]
\label{thm-EH}
Let $(M,g)$ be a complete noncompact Riemannian manifold which satisfies the uniform volume growth condition
\[
\vol_t(B_r(p)) \le e^{k(1+r^2)}
\]
for some point $p\in M$ and a uniform constant $k > 0$ for all $t\in [0,T]$. Let $w \in \calC^\infty(M\times (0,T]) \cap
\calC^0(M\times [0,T])$ satisy the differential inequality
\[
\frac{\partial}{\partial t} w \le \Delta w + {\bf a}\cdot\nabla w + b w,
\]
where $\sup_{M\times[0,T]}|{\bf a}| \le c_1$ and $\sup_{M\times[0,T]}|b| \le c_2$. If in addition
$w(x,0) \le 0$ for all $x\in M$, $\int_0^T\int_M e^{-c_3 r_t(p,y)^2}|\nabla w|^2(y)\, d\mu_t\, dt < \infty$ for some constant $c_3 > 0$
and $\sup_{M\times[0,T]}|\frac{\partial}{\partial t} g(t)| \le c_4$, then $w \le 0$ on $M\times [0,T]$.
\end{thm} 

We  now state a Harnack estimate which holds for $\tilde{R}$, allowing one to compare its value at different space-time points as the metric evolves under Ricci flow. 
Hamilton first proved the Harnack estimate for Ricci flow in the case that the curvature operator is nonnegative. Chow \cite{Ch} 
generalized this to the case which allows some negative curvature. Using the maximum principles stated above, and using the 
uniform boundedness in space  of $\tilde{R}$ and its derivatives on a fixed time slice, one can adapt Chow's arguments to the 
case of Ricci flow on complete manifolds. 
The Harnack estimate for $\tilde{R}$ states that there exist constants $\tau_0 > 0$ and $C$ such that for every $x_1, x_2 \in M$ and $\tau_2 \ge \tau_1 \ge \tau_0$,
\begin{equation}
\label{eq-har-r}
\tilde{R}(x_2,\tau_2) + 1 \ge e^{-\Delta/4 - C(\tau_2 - \tau_1)}(\tilde{R}(x_1,\tau_1) + 1),
\end{equation}
where 
$$\Delta = \Delta(x_1,x_2,\tau_1,\tau_2) = \inf_{\gamma}\int_{\tau_1}^{\tau_2} |\frac{d\gamma}{dt}(t)|^2\, dt,$$
and the infimum is taken over all paths $\gamma$ in $M$ whose graphs $(\gamma(t),t)$ join $(x_1,\tau_1)$ and $(x_2,\tau_2)$. Since the metric $\tilde{g}$ is shrinking it is easy to derive (see \cite{Ha2}) 
\begin{equation}
\label{eq-delta}
\Delta  \le \frac{d(x_1,x_2,\tau_1)^2}{\tau_2-\tau_1} \le \frac{\dist(x_1,x_2,0)^2}{\tau_2-\tau_1},
\end{equation}
where $d(x_1,x_2,\tau_1)$ is a distance between points $x_1$ and $x_2$ computed at time $\tau_1$.

\end{document}